\documentclass[a4paper,11pt]{article}
\usepackage{tipa}
\usepackage{bbm}
\usepackage{bbding}
\usepackage{stmaryrd}
\usepackage{amsfonts}
\usepackage{graphicx}
\usepackage{epsfig}
\usepackage{amssymb}
\usepackage{amssymb,amsmath,amsthm}
\usepackage{pifont}
\usepackage{graphicx}
\usepackage{indentfirst}
\theoremstyle{definition}

\numberwithin{equation}{section}

\DeclareSymbolFontAlphabet{\testi}{letters} \large\normalsize
\def\dis{\displaystyle}

\begin{document}

\title{\bf Delta-shocks and vacuums in zero-pressure gas dynamics by the flux approximation
\thanks{\ Supported by the NSF of China (11361073).} }

 \author{Hanchun Yang \ \ Jinjing Liu \ \\\
 \footnotesize\slshape{Department of Mathematics, Yunnan University, Kunming 650091, P.R. China}}

\date{}\maketitle

{\small

\noindent{\bf Abstract:} In this paper, firstly, by solving the
Riemann problem of the zero-pressure flow in gas dynamics with a
flux approximation, we construct parameterized delta-shock and
constant density solutions, then we show that, as the flux
perturbation vanishes, they converge to the delta-shock and vacuum
state solutions of the zero-pressure flow, respectively. Secondly,
we solve the Riemann problem of the Euler equations of isentropic
gas dynamics with a double parameter flux approximation including
pressure. Further we rigorously prove that, as the two-parameter
flux perturbation vanishes, any Riemann solution containing two
shock waves tends to a delta shock solution to the zero-pressure
flow; any Riemann solution containing two rarefaction waves tends to
a two-contact-discontinuity solution to the zero-pressure flow
and the nonvacuum intermediate state in between tends to a vacuum
state.

\vspace{0.1cm}

\noindent{\it Keywords}:  Euler equations of isentropic gas
dynamics; Zero-pressure flow; Transport equations; Riemann problem;
Delta shock wave; Vacuum; flux approximation.

\vspace{0.1cm}

\noindent{\it AMS subject classifications}: 35L65, 35B30, 76E19, 35Q35, 35L67

}

\date{}

\section{ Introduction }

The well-known zero-pressure gas dynamics reads
\begin{align}\label{eq1.1}
\left\{\begin{array}{l}
  \rho_t  + (\rho u)_x  =0,\\
  (\rho u)_t  +  (\rho u^2)_x =0,
\end{array}\right.
\end{align}
which are also called the transport equations, or Euler equations
for pressureless fluids, where $\rho$ is the density and $u$ the
velocity. It can be used to model the motion of free
particles which stick under collision \cite{WE,BGre} and the
formation of large-scale structures in the universe \cite{S-Z}.

In the past twenty years, there has been a great explosion of
interests in the extensive investigations on the zero-pressure gas
dynamics, for instance,
see \cite{Bouchut,WE,BGre,Sheng-Zhang,Li-Zhang2,Li-Yang,Huang-Wang},
etc. Among these works, Bouchut \cite{Bouchut} first established the
existence of measure solutions of the Riemann problem. Weinan E,
Rykov and Sinai \cite{WE} studied the existence of global weak
solution and the behavior of such global solution with random
initial data. The 1-D and 2-D Riemann problems were solved by Sheng
and Zhang \cite{Sheng-Zhang} with the characteristic analysis and
the vanishing viscosity method, see also \cite{Li-Zhang2}.
Huang and Wang \cite{Huang-Wang} obtained the uniqueness result of
weak solution when the initial data is a Radon measure. In these
papers it has been proved that $\delta$-shock waves and vacuum
states do occur in solutions. For $\delta$-shock waves, we refer to
\cite{Ko,Keyfitz,Floch,Tan-Zhang,Tan-Zhang-Zheng,Yang1,Yang-Zhang1,Yang-Zhang2}  for more details.

During the recent decade, the problem concerning the phenomena of concentration and
cavitation and the formation of $\delta$-shock waves and vacuum
states in solutions has received much attention. For example, see
\cite{C-L-1,C-L-2,Ljq,Yin-Sheng,Shen-Sun,Yang-Wang,Cheng-Yang}, etc.
In 2003, Chen and Liu \cite{C-L-1} considered the Euler
equations of isentropic gas dynamics
\begin{align}\label{eq1.2}
\left\{\begin{array}{l}
   \rho_t  + (\rho u)_x  =0, \cr\noalign {\vskip1truemm}
(\rho u)_t  +  (\rho u^2 + P)_x =0,
\end{array}\right.
\end{align}
where $\rho\geq0$, $u$, $P$ denote the density, the velocity and the
pressure respectively. The scalar pressure $P(\rho, \epsilon)$
satisfies
\begin{align}
\lim_{\epsilon\rightarrow 0}P(\rho, \epsilon)=0,
\end{align}
where $\epsilon>0$ is a small parameter. In their works, in
\eqref{eq1.2} Chen and Liu took the prototypical pressure functions
for ploytropic gas
\begin{align}
P(\rho, \epsilon)=\epsilon p(\rho),\ \ p(\rho)=\rho^\gamma/\gamma,\
\ \ \gamma>1.
\end{align}
They identified and analyzed the phenomena of concentration and
cavitation and the formation of $\delta$-shock waves and vacuum
states in solutions to the system \eqref{eq1.2} with (1.4) as
$\epsilon\rightarrow 0$. Further, in \cite{C-L-2} they also studied
the nonisentropic fluids. Specially, Li \cite{Ljq} investigated
the zero temperature limit for $\gamma=1$ in (1.4).
Besides, the results were extended to the relativistic
Euler equations for polytropic gases by Yin and Sheng
\cite{Yin-Sheng}, the perturbed Aw-Rascle model by Shen and Sun
\cite{Shen-Sun}, etc. Very recently, see \cite{Yang-Wang,Cheng-Yang}
for the modified Chaplygin gas pressure law.  All in all, these
works on this topic are only focused on the pressure level.

Motivated partly by \cite{C-L-1,C-L-2,Lu}, in the present paper, by
introducing a flux approximation, we propose to consider the
following system
\begin{align}\label{eq1.5}
\left\{\begin{array}{l}
 \rho_t+(\rho u-2\epsilon_1
u)_x=0,\cr\noalign{\vskip2truemm}
 (\rho u)_t+\big(\rho u^2-\epsilon_1u^2+\epsilon_2 p(\rho)\big)_x=0,
\end{array}\right.
\end{align}
where the density $\rho\geq2\epsilon_1$, $\epsilon_1,\epsilon_2 >0$
are parameters. Physically, a reasonable perturbation can be used to
govern some dynamical behaviors of fluids, so it is worth studying
the flux perturbation problem which plays an important role in all
the three of theory, application and computation. In contrast to the
previous works in \cite{C-L-1,C-L-2,Ljq,Yin-Sheng,Shen-Sun,Yang-Wang,Cheng-Yang}, we here develop a flux approximation
approach which contains the pressure perturbation portion.

Firstly we consider a special case $\epsilon_2 =0$ in \eqref{eq1.5},
that is
\begin{align}\label{eq1.6}
\left\{\begin{array}{l}
 \rho_t+(\rho u-2\epsilon_1
u)_x=0,\cr\noalign{\vskip2truemm}
 (\rho u)_t+\big(\rho u^2-\epsilon_1u^2\big)_x=0,
\end{array}\right.
\end{align}
this is a pure flux approximation of special curiosity. We solve the
Riemann problem of the system \eqref{eq1.6} with initial conditions
\begin{align}\label{eq1.7}
(\rho,u)(0,x)=\left\{\begin{array}{ll}
  (\rho_-,u_-),&x<0,\cr\noalign{\vskip1truemm}
  (\rho_+,u_+), &x>0,
 \end{array}\right.
\end{align}
where $(\rho_\pm,u_\pm)$ are arbitrary constants. The Riemann
solutions include two kinds of somewhat interesting features. When
$u_-<u_+$, the solution consists of two contact discontinuities and
a constant density state besides two constant states. When
$u_->u_+$, the solution contains a delta shock wave depending on a
parameter. From the solutions constructed, one can find that,
compared with the zero-pressure gas dynamics, the vacuum state is
removed, while for the $\delta$-shock wave, the location and
propagation speed are preserved, the weight decreases. Theses mean
that the flux perturbation works in the pressureless gases.

Then we prove that, as the flux approximation vanishes, that is,
parameter $\epsilon_1\rightarrow0$, any parameterized delta-shock
solution converges to the corresponding one of the zero-pressure
 flow \eqref{eq1.1}. By contrast, any constant density solution goes to the vacuum
solution.

Secondly, we solve the Riemann problem \eqref{eq1.5},
\eqref{eq1.7}. Because both of the characteristic fields are
genuinely nonlinear, the elementary waves consist of backward
centred rarefaction wave ($\overleftarrow{R}$), forward centred
rarefaction wave ($\overrightarrow{R}$), backward shock wave
($\overleftarrow{S}$) and forward shock wave ($\overrightarrow{S}$).
The curves of elementary waves divide the phase plane into five
domains. By the analysis method in phase plane, we can establish the
existence and uniqueness of Riemann solutions including five
different structures.

Moreover, we analyze the limit of Riemann solutions of (1.5) and
(1.7) as the double parameter ${\epsilon_1,\epsilon_2\rightarrow0}$.
It is shown that when $u_+< u_-$, the Riemann solution containing
two shock waves converges to a delta shock solution, which is
exactly the solution to zero-pressure  flow \eqref{eq1.1}. The density between the
two shock waves tends to an extreme concentration in the form of a
weighted $\delta$-function, which results in the formation of a
delta shock wave. Besides, it is also shown that when $u_+> u_-$,
the Riemann solution containing two rarefaction waves tends to a
two-contact-discontinuity solution to zero-pressure  flow \eqref{eq1.1}, and the
nonvacuum intermediate state in between tends to a vacuum state as
${\epsilon_1,\epsilon_2\rightarrow0}$.

Following the above analysis, one can find a fact of interest, that is, the flux approximations of difference have their respective effect on the formation of delta-shock and vacuum state in isentropic fluids. In this regard, it is different from those only in pressure level \cite{C-L-1,C-L-2,Ljq,Yin-Sheng,Shen-Sun,Yang-Wang,Cheng-Yang}. Meanwhile, the results obtained show that both the delta shock wave and vacuum are stable under some flux small perturbations. Therefore this work extends in some sense the previous results and proofs in \cite{C-L-1,C-L-2,Ljq}.
The flux approximation method can be also extend to the Euler equations for nonisentropic fluids and Chaplygin gas equations \cite{Yang-Liu,Liu-Yang}.

The arrangement of this paper is as follows. In Section 2, we recall
the solutions of (1.1), (1.7). Section 3 solves the Riemann problem
(1.6), (1.7) and discusses the limits of Riemann solutions. Section
4 solves the Riemann problem for \eqref{eq1.5}. Sections
5 and 6 investigate the limit of solutions of (1.5) and (1.7).

\section{ Delta-shocks and vacuums for the zero-pressure flow}

As a start, we briefly recall $\delta$-shocks and vacuum states in
the Riemann solutions to the zero-pressure flow \eqref{eq1.1}. See
\cite{Sheng-Zhang,Li-Zhang2} for more details.

The system \eqref{eq1.1} has a double eigenvalue $\lambda=u$ with
the associated eigenvector $r=(1,0)^T$ satisfying
$\nabla\lambda\cdot r=0$, which means that it is nonstrictly
hyperbolic and $\lambda$ linearly degenerate.

Consider Riemann problem \eqref{eq1.1},\eqref{eq1.7}. By seeking
self-similar solution $(\rho,u)(t,x)=(\rho,u)(\xi)\ (\xi=x/t)$, it
is easy to find that, besides the constant state and singular
solution $\rho=0, u=\xi$(vacuum state), the elementary waves
of \eqref{eq1.1} are nothing but contact discontinuities. The
Riemann problem can be solved by the
following two cases.

For the case $u_- < u_+$, the solution includes two contact
discontinuities and a vacuum state besides constant states. That is,
\begin{align}\label{eq2.1}
 (\rho,u)(\xi)=\left\{\begin{array}{lc}
  (\rho_-,u_-),&-\infty<\xi< u_-,\\
   (0,\xi),&u_-\leq\xi\leq u_+,\\
   (\rho_+,u_+),&u_+<\xi<+\infty.
   \end{array}\right.
\end{align}

For the case $u_- > u_+$, a solution containing a weighted
$\delta$-measure (i.e., $\delta$-shock) supported on a line will
develop in solutions due to the overlap of characteristic lines.

To define the measure solution, a two-dimensional weighted
$\delta$-function $w(s)$$\delta_S$ supported on a smooth curve $S$
parameterized as $t=t(s)$, $x=x(s) (c\leq s\leq d)$ can be defined
by
\begin{align}\label{eq2.2}
\langle w(t(s))\delta_S, \varphi(t(s),x(s))\rangle =\dis\int^d_c\dis
w(t(s))\varphi(t(s),x(s))\sqrt{x'(s)^2 + t'(s)^2}ds
\end{align}
for all the test functions $\varphi(t,x)\in C^\infty_0 (R^+ \times
R^1)$.

With this definition, a $\delta$-shock solution of \eqref{eq1.1} can
be introduced as follows
\begin{align}\label{eq2.3}
\begin{array}{l}
\rho(t,x)=\rho_0(t,x)+w(t)\delta_S, \ \ \ \ u(t,x)=u_0(t,x),
\end{array}
\end{align}
where $S=\{(t,\sigma t):0\leq t <\infty\}$,
\begin{align}\label{eq2.4}
\begin{array}{ll}
\rho_0(t,x)=\rho_-+[\rho]\chi(x-\sigma t),
u_0(t,x)=u_-+[u]\chi(x-\sigma t),
w(t)=\dis\frac{t}{\sqrt{1+\sigma^2}}(\sigma[\rho]-[\rho u]),
\end{array}
\end{align}
in which $[g]=g_+-g_-$, $\sigma$ is the velocity of the $\delta$-shock, and
$\chi(x)$ the characteristic function that is 0 when $x<0$ and 1
when $x>0$.

As shown in \cite{Sheng-Zhang,Li-Zhang2}, for any $\varphi(t,x)\in C^\infty_0
(R^+ \times R^1)$, the $\delta$-shock solution constructed above
satisfies
\begin{align}\label{eq2.5}
\left\{\begin{array}{l}
 \langle \rho, \varphi_t\rangle +\langle \rho
u,\varphi_x\rangle=0,\cr\noalign{\vskip2truemm} \langle \rho
u,\varphi_t\rangle+\langle \rho u^2,
\varphi_x\rangle=0,\cr\noalign{\vskip2truemm}
\end{array}\right.
\end{align}
where
\begin{align}\label{eq2.6}
\begin{array}{l}
 \langle \rho,
\varphi\rangle=\dis\int^{+\infty}_0\dis\int^{+\infty}_{-\infty}\dis\rho_0
\varphi dx dt+\langle
w\delta_S,\varphi\rangle,\cr\noalign{\vskip2truemm} \langle\rho
u,\varphi\rangle=\dis\int^{+\infty}_0\dis\int^{+\infty}_{-\infty}\dis\rho_0
u_0 \varphi dx dt+\langle \sigma
w\delta_S,\varphi\rangle.\cr\noalign{\vskip2truemm}
\end{array}
\end{align}

Furthermore, substituting \eqref{eq2.3} and \eqref{eq2.4} into
\eqref{eq2.5} under the condition \eqref{eq2.2} and \eqref{eq2.6},
one can get the generalized Rankine-Hugoniot relation
\begin{equation}\label{eq2.7}
\dis\frac{dx}{dt}=\sigma,\
\dis\frac{d\big(w(t)\sqrt{1+\sigma^2}\big)}{dt}=\sigma[\rho]-[\rho
u], \
\dis\frac{d\big(w(t)\sigma\sqrt{1+\sigma^2}\big)}{dt}=\sigma[\rho
u]-[\rho u^2]
\end{equation}
which reflects the relationship among the
location, weight and propagation speed of the $\delta$-shock wave.

To guarantee the uniqueness, the entropy condition is supplemented
as
\begin{align}\label{eq2.8}
u_+<\sigma< u_-,
\end{align}
which means that all characteristic lines on both sides of the
discontinuity are not out-going. So it is a overcompressive
condition.

Then solving the generalized Rankine-Hugoniot relation \eqref{eq2.7}
with initial data $x(0)=0$ and $w(0)=0$ under the entropy condition
\eqref{eq2.8} yields
\begin{align}\label{eq2.9}
\begin{array}{l}
\sigma=\dis\frac{\sqrt{\rho_+}u_+ +
\sqrt{\rho_-}u_-}{\sqrt{\rho_+}+\sqrt{\rho_-}}, \
w(t)=\dis\frac{\sqrt{\rho_+ \rho_-}(u_- - u_+)t}{\sqrt{1+\sigma
^2}}.
\end{array}
\end{align}
Therefore, a $\delta$-shock solution defined by \eqref{eq2.3} with
\eqref{eq2.4} and \eqref{eq2.9} is obtained.

\section {\bf Riemann solutions and limit analysis of \eqref{eq1.6} as $\epsilon_1\rightarrow0$ }

The section solves the Riemann problem \eqref{eq1.6}, \eqref{eq1.7},
and studies the limit of solutions.

For the system \eqref{eq1.6}, the eigenvalue and the associated
eigenvector are $\lambda =u$ and $r=(1,0)^T $, respectively,
satisfying $\nabla\lambda \cdot r=0$, which means that the system
\eqref{eq1.6} is full linear degenerate and elementary waves only
involve contact discontinuities.

In a similar way as the Riemann problem (1.1), (1.7), it is easy to
find that, for smooth solutions, besides the constant state, the
system \eqref{eq1.6} provides the singular solution
\begin{align}\label{eq3.3}
 \rho=2\epsilon_1,\ u= \xi,
\end{align}
which is called constant density states. While the elementary wave
has only contact discontinuity
\begin{align}\label{eq3.5}
J:\ \ \omega=\xi=u_-=u_+,
\end{align}
which is characterized by $x/t=u_-=u_+$ in $(t,x)$-plane. It can
connect two states $(\rho_-,u_-)$ and $(\rho_+,u_+)$ if and only if
they are located on the line $u=u_-=u_+$ in the $(\rho,u)$-plane.

Now, with constants, constant density state and contact
discontinuity, we construct the solutions of Riemann problem \eqref{eq1.6},
\eqref{eq1.7} by two cases.

For the case $u_-<u_+$, we draw lines $u=u_-$ and $u=u_+$ from
$(\rho_-,u_-)$ and $(\rho_+,u_+)$, respectively, in the
$(\rho,u)$-plane. These two lines intersect the line
$\rho=2\epsilon_1$ at $(2\epsilon_1,u_-)$ and $(2\epsilon_1,u_+)$.
Thus the solution can be constructed by two contact discontinuities
and a constant-density state besides two constant states (see Fig.
1), and can be expressed as
\begin{align}\label{eq3.6}
 (\rho,u)(t,x)=(\rho,u)(\xi)=\left\{\begin{array}{lc}
  (\rho_-,u_-),&-\infty<\xi< u_-,\\
   (2\epsilon_1,\xi),&u_-\leq\xi\leq u_+,\\
   (\rho_+,u_+),&u_+<\xi<+\infty.
   \end{array}\right.
\end{align}

\begin{center}
\includegraphics[width=2.0in]{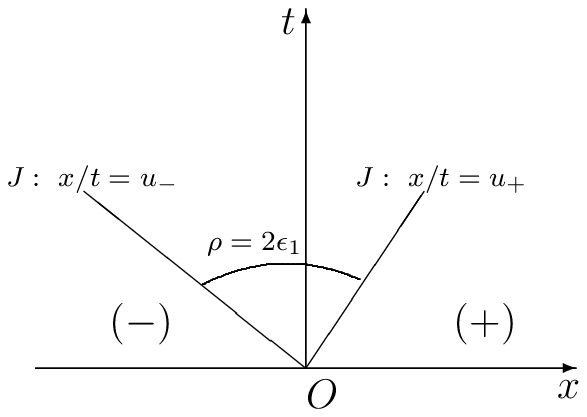}\hspace{2.0cm}\includegraphics[width=1.4in]{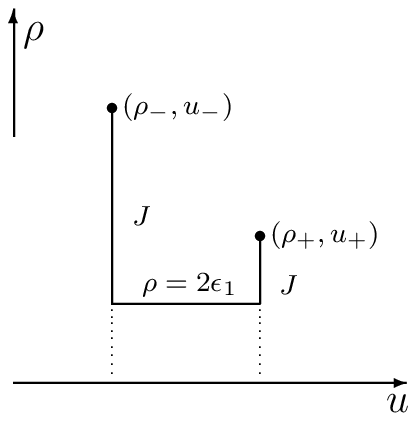}

\vspace{-0.15cm}

\hspace*{0.5cm}{\small Fig. 1. Constant density.}

\end{center}

\begin{center}
\includegraphics[width=2.0in]{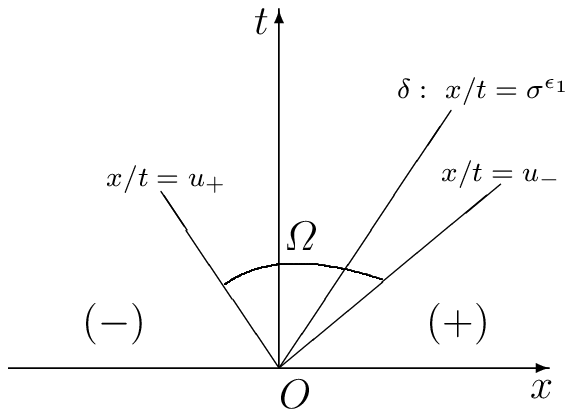}
\vspace{-0.15cm}

\hspace*{0.5cm}{\small Fig. 2. Characteristic analysis of
$\delta$-shock.}

\end{center}

For the case $u_- > u_+$, as indicated in Fig. 2, since characteristic lines from initial data
overlap each other in the region $\Omega$, so the singularity of solutions must develop in this region.
As shown in \cite{Sheng-Zhang,Li-Zhang2,Yang1}, there is no solutions
 exist in bounded variation space, then we can construct the Riemann solution by a delta-shock wave.

With the definitions as the above section, we seek a delta shock
solution $(\rho^{\epsilon_1}, u^{\epsilon_1}, \sigma^{\epsilon_1},
w^{\epsilon_1})$ of the form (2.3), (2.4), then the following
generalized Rankine-Hugoniot relation holds,
\begin{align}\label{eq3.17}
\left\{\begin{array}{l}
\dis\frac{dx}{dt}=\sigma^{\epsilon_1},\cr\noalign{\vskip2truemm}
\dis\frac{d\big(w^{\epsilon_1}(t)\sqrt{1+(\sigma^{\epsilon_1})^2}\big)}{dt}=\sigma^{\epsilon_1}[\rho]-[\rho
u-2\epsilon_1 u],\cr\noalign{\vskip2truemm}
\dis\frac{d\big(w^{\epsilon_1}(t)\sigma^{\epsilon_1}
\sqrt{1+(\sigma^{\epsilon_1})^2}\big)}{dt}=\sigma^{\epsilon_1}[\rho
u]-[\rho u^2-\epsilon_1 u^2],
\end{array}\right.
\end{align}

Besides, the discontinuity should satisfy the entropy condition
\begin{equation}\label{eq3.19}
u_+<\sigma^{\epsilon_1}<u_-.
\end{equation}

In what follows, the generalized Rankine-Hugoniot relation will be
applied in particular to Riemann problem \eqref{eq1.6} and
\eqref{eq1.7} for the case $u_->u_+$. Now this Riemann problem is
reduced to solving \eqref{eq3.17} with the initial conditions $t=0:\
x(0)=0,\ w^{\epsilon_1}(0)=0.$

Obviously, we have from \eqref{eq3.17} that
\begin{align}\label{eq3.21}
\left\{\begin{array}{l}
w^{\epsilon_1}(t)\sqrt{1+(\sigma^{\epsilon_1})^2}=[\rho]x-[\rho
u-2\epsilon_1 u]t,\cr\noalign{\vskip2truemm}
w^{\epsilon_1}(t)\sigma^{\epsilon_1}\sqrt{1+(\sigma^{\epsilon_1})^2}=[\rho
u]x-[\rho u^2-\epsilon_1 u^2]t.
\end{array}\right.
\end{align}
Multiplying the first equation by $\sigma^{\epsilon_1}$ and together
with the second equation to give
\begin{equation}\label{eq3.23}
\dis\frac{d([\rho]\frac{x^2}{2}-[\rho u]xt)}{dt}=-2[\epsilon_1
u]\sigma^{\epsilon_1} t-[\rho u^2-\epsilon_1 u^2]t.
\end{equation}
In view of the knowledge concerning delta shock waves in
\cite{Sheng-Zhang,Li-Zhang2,Yang1}, we find that
$\sigma^{\epsilon_1}$ is a constant. Then it follows from
\eqref{eq3.23} that
\begin{equation}\label{eq3.24}
\frac{[\rho]}{2}x^2-[\rho u]tx+([2\epsilon_1
u]\sigma^{\epsilon_1}+[\rho u^2-\epsilon_1 u^2])\frac{t^2}{2}=0,
\end{equation}
one solves
\begin{equation}\label{eq3.25}
x(t)=\dis\frac{[\rho u]\pm\sqrt{[\rho u]^2-[\rho]([2\epsilon_1
u]\sigma^{\epsilon_1}+[\rho u^2-\epsilon_1 u^2])}}{[\rho]}t,
\end{equation}
as $[\rho]\neq0$. Then we can obtain
\begin{equation}\label{eq3.28}
\sigma^{\epsilon_1}=\dis\frac{[(\rho-\epsilon_1)
u]+\sqrt{(\rho_--\epsilon_1)(\rho_+-\epsilon_1)}(u_--u_+)}{[\rho]},
\end{equation}
under the entropy condition \eqref{eq3.19}. Therefore, from (3.6) we
get
\begin{equation}\label{eq3.29}
w^{\epsilon_1}(t)=\dis\frac{[\epsilon_1
u]+\sqrt{(\rho_--\epsilon_1)(\rho_+-\epsilon_1)}(u_--u_+)}{\sqrt{1+(\sigma^{\epsilon_1})^2}}t.
\end{equation}
Especially, when $[\rho]=0$,
\begin{equation}\label{eq3.30}
\sigma^{\epsilon_1}=\dis\frac{u_-+u_+}{2},\
x(t)=\dis\frac{u_-+u_+}{2}t,\
w^{\epsilon_1}(t)=\dis\frac{[2\epsilon_1 u-\rho
u]}{\sqrt{1+(\sigma^{\epsilon_1})^2}}t.
\end{equation}

With the above analysis, we reach the following result.
\vspace{0.2cm}

\noindent{\bf Theorem 3.1.} \textit{The Riemann problem
\eqref{eq1.6},\eqref{eq1.7} admits a unique weak solution which
includes a constant density state as $u_-<u_+$ and a
$\delta$-shock wave as $u_->u_+$.}

\vspace{0.2cm}

In the next, we proceed to discuss the limit of Riemann solutions of
the system \eqref{eq1.6} as $\epsilon_1\rightarrow0$ for $\rho_-\neq
\rho_+$. It needs to investigate two cases: $(1)\ u_->u_+$, $(2)\
u_-<u_+$.

We first consider Case (1), which is relevant to the formation of
$\delta$-shock waves.

Computing the limits of $\sigma^{\epsilon_1}$ and
$w^{\epsilon_1}(t)$ as $\epsilon_1\rightarrow0$, from \eqref{eq3.28}
and \eqref{eq3.29}, one can obtain
\begin{equation}\label{eq4.1}
\lim_{\epsilon_1\rightarrow
0}\sigma^{\epsilon_1}=\dis\frac{\sqrt{\rho_+}u_+ +
\sqrt{\rho_-}u_-}{\sqrt{\rho_+}+\sqrt{\rho_-}}=\sigma,\
\lim_{\epsilon_1\rightarrow 0}w^{\epsilon_1}(t)
 =\dis\frac{\sqrt{\rho_+
\rho_-}(u_--u_+)}{\sqrt{1+\sigma^2}}t=w(t).
\end{equation}

In a simple way similar to that in \cite{C-L-1,Ljq}, one can
easily conclude the following results.

\vspace{0.2cm}

\noindent{\bf Theorem 3.2.} \textit{Let $u_->u_+$. For each fixed
$\epsilon_1>0$, assume that $(\rho^{\epsilon_1},u^{\epsilon_1})$ is
a $\delta$-shock solution of \eqref{eq1.6}, \eqref{eq1.7}. Then,
when $\epsilon_1\rightarrow 0$, the pair of limit functions
$(\rho,u)$ is a $\delta$-shock solution of \eqref{eq1.1},
\eqref{eq1.7}. Moreover, $\rho$ and $\rho u$ are the sum of a step
function and a $\delta$-measure with weights
$\dfrac{t}{\sqrt{1+\sigma^2}}(\sigma[\rho]-[\rho u])$ and
$\dfrac{t}{\sqrt{1+\sigma^2}}(\sigma[\rho u]-[\rho u^2])$,
respectively.}

\vspace{0.2cm}

%\noindent{\bf Remark 3.3.} According to the  Theorem 3.2, it is easy
%to see that the generalized Rankine-Hugoniot relation \eqref{eq2.7}
%holds.

Now we turn to Case (2). At this moment, the solution of Riemann
problem \eqref{eq1.6}, \eqref{eq1.7} can be expressed as
\eqref{eq3.6}. It is obvious to get that, as
$\epsilon_1\rightarrow0$, the limit of solution to the system
\eqref{eq1.6} is just the vacuum solution \eqref{eq2.1} to the
zero-pressure flow.

\section {\bf Solutions of Riemann problem \eqref{eq1.5} and \eqref{eq1.7}}

In this section, we solve the elementary waves and  construct the
solutions of Riemann problem \eqref{eq1.5} and \eqref{eq1.7}.
For any $\epsilon_1,\epsilon_2>0$, the system \eqref{eq1.5} has two
eigenvalues
\begin{align}\label{eq4.1}
\begin{array}{l}
\lambda_1=u-\sqrt{\epsilon_2 \rho^{\gamma-2} (\rho-2\epsilon_1)}\ ,
\qquad  \lambda_2=u+\sqrt{\epsilon_2 \rho^{\gamma-2}
(\rho-2\epsilon_1)},
\end{array}
\end{align}
so it is strictly hyperbolic. The corresponding right eigenvectors
are
$$
\begin{array}{l}
{r}_1=\Bigg(1, \dis-\sqrt{\frac{\epsilon_2
\rho^{\gamma-2}}{\rho-2\epsilon_1}} \Bigg)^T, \qquad {r}_2=\Bigg(1,
\dis\sqrt{\frac{\epsilon_2 \rho^{\gamma-2}}{\rho-2\epsilon_1}}
\Bigg)^T.
\end{array}
$$
Since $\nabla\lambda_i\cdot {r}_i\neq0\ (i=1,2)$, both of the
characteristic fields are genuinely nonlinear.

Seeking the self-similar solution, we reach the following boundary
value problem
\begin{align}\label{eq4.2}
\left\{\begin{array}{l}
 -\xi\rho_\xi+(\rho u-2\epsilon_1 u)_\xi=0,\cr\noalign{\vskip3truemm}
 -\xi(\rho u)_\xi+(\rho u^2-\epsilon_1 u^2+\dis\frac{\epsilon_2\rho^{\gamma}}{\gamma})_\xi=0,
\end{array}\right.
\end{align}
and
\begin{align}\label{eq4.3}
(\rho,u)(\pm\infty)=(\rho_\pm,u_\pm).
\end{align}

For any smooth solution, \eqref{eq4.2} is equivalent to
\begin{align}\label{eq4.4}
\left(\begin{array}{cc}
 -\xi+u &\rho-2\epsilon_1\cr\noalign
{\vskip2truemm}
 -\xi u+u^2+\epsilon_2\rho^{\gamma-1} & -\xi \rho+2\rho
 u-2\epsilon_1 u
\end{array}\right)
\left(\begin{array}{c}
 d\rho\cr\noalign {\vskip2truemm}
  du
\end{array}\right)=0,
\end{align}
which provides either the general solution (constant state)
\begin{align}\label{eq4.5}
(\rho,u)(\xi)=constant,
\end{align}
or the backward centred rarefaction wave
\begin{align}\label{eq4.6}
\overleftarrow{R}(\rho_-,u_-): \ \ \left\{\begin{array}{l}
 \xi=\lambda_1=u-\sqrt{\epsilon_2 \rho^{\gamma-2} (\rho-2\epsilon_1)},\cr\noalign {\vskip2truemm}
 u-u_-=-\dis\int^{\rho}_{\rho_-} \sqrt{\frac{\epsilon_2
s^{\gamma-2}}{s-2\epsilon_1}} ds,\hspace{1cm} \rho<\rho_-,
\end{array}\right.
\end{align}
or the forward centred rarefaction wave
\begin{align}\label{eq4.7}
\overrightarrow{R}(\rho_-,u_-):\ \ \left\{\begin{array}{l}
 \xi=\lambda_2=u+\sqrt{\epsilon_2 \rho^{\gamma-2} (\rho-2\epsilon_1)},\cr\noalign {\vskip2truemm}
 u-u_-=\dis\int^{\rho}_{\rho_-} \sqrt{\frac{\epsilon_2
s^{\gamma-2}}{s-2\epsilon_1}} ds,\hspace{1.4cm} \rho>\rho_-.
\end{array}\right.
\end{align}

For the backward centred rarefaction wave, differentiating
$u$ with respect to $\rho$ in the second equation of \eqref{eq4.6},
it follows that $\dis u_{\rho}=-\sqrt{\frac{\epsilon_2
\rho^{\gamma-2}}{\rho-2\epsilon_1}}<0$. For the forward centred
rarefaction wave, it is easy to see that $\dis
u_{\rho}=\sqrt{\frac{\epsilon_2
\rho^{\gamma-2}}{\rho-2\epsilon_1}}>0$.

Taking the limit $\rho\rightarrow 2\epsilon_1$ in the second
equation of \eqref{eq4.6} leads to
\begin{align}\label{eq4.8}
\lim_{\rho\rightarrow 2\epsilon_1}u=u_-+\dis
\int^{\rho_-}_{2\epsilon_1} \sqrt{\frac{\epsilon_2
s^{\gamma-2}}{s-2\epsilon_1}} ds.
\end{align}
Since $\underset{s\rightarrow 2\epsilon_1}\lim
\Bigg((s-2\epsilon_1)^{\frac{1}{2}}\dis\sqrt{\frac{\epsilon_2
s^{\gamma-2}}{s-2\epsilon_1}}\Bigg)=\sqrt{\epsilon_2(2\epsilon_1)^{\gamma-2}}$,
the integral $\dis \int^{\rho_-}_{2\epsilon_1}
\sqrt{\frac{\epsilon_2 s^{\gamma-2}}{s-2\epsilon_1}} ds$ is
convergent due to Cauchy criterion. Thus, from \eqref{eq4.8}, we can conclude that
the backward centred rarefaction wave curve intersects with the line
$\rho=2\epsilon_1$ at the point
$(2\epsilon_1,u_1)=(2\epsilon_1,u_-+\dis \int^{\rho_-}_{2\epsilon_1}
\sqrt{\frac{\epsilon_2 s^{\gamma-2}}{s-2\epsilon_1}} ds)$.

Performing the limit $\rho\rightarrow +\infty$ in the second
equation in \eqref{eq4.7} yields
\begin{align}\label{eq4.9}
\lim_{\rho\rightarrow +\infty}u=u_-+\dis \int^{+\infty}_{\rho_-}
\sqrt{\frac{\epsilon_2 s^{\gamma-2}}{s-2\epsilon_1}} ds.
\end{align}
Since $\sqrt{\frac{\epsilon_2 s^{\gamma-2}}{s-2\epsilon_1}}>\sqrt{\frac{\epsilon_2 s^{\gamma-2}}{s}}$, we have
\begin{align}\label{eq4.10}
\dis \int^{+\infty}_{\rho_-} \sqrt{\frac{\epsilon_2
s^{\gamma-2}}{s-2\epsilon_1}} ds>\dis \int^{+\infty}_{\rho_-}
\sqrt{\frac{\epsilon_2 s^{\gamma-2}}{s}} ds=+\infty.
\end{align}
Thus, from \eqref{eq4.9}, one deduces that $\underset{\rho\rightarrow
+\infty}\lim u=+\infty$.

For a bounded discontinuity at  $\xi=\sigma^{\epsilon_1\epsilon_2}$,
the Rankine-Hugoniot relation
\begin{align}\label{eq4.11}
\left\{\begin{array}{l}
 -\sigma^{\epsilon_1\epsilon_2}[\rho]+[\rho u-2\epsilon_1 u]=0,\cr\noalign {\vskip3truemm}
 -\sigma^{\epsilon_1\epsilon_2}[\rho u]+[\rho u^2-\epsilon_1 u^2+\dis\frac{\epsilon_2\rho^{\gamma}}{\gamma}]=0,
\end{array}\right.
\end{align}
holds, where $[q]=q_r-q_l$ with $q_l=q(t,x(t)-0)$ and
$q_r=q(t,x(t)+0)$.

Eliminating $\sigma^{\epsilon_1\epsilon_2}$ from \eqref{eq4.11}, we
get
\begin{align}\label{eq4.12}
\big(\rho_l\rho_r-\epsilon_1(\rho_l+\rho_r)\big)(u_r-u_l)^2=\frac{\epsilon_2}{\gamma}(\rho_r-\rho_l)(\rho_r^{\gamma}-\rho_l^{\gamma}),
\end{align}
which yields $\rho_l\rho_r-\epsilon_1(\rho_l+\rho_r)>0$. Thus, we have
\begin{align}\label{eq4.13}
u_r-u_l=\dis\pm\sqrt{\frac{\epsilon_2(\rho_r-\rho_l)(\rho_r^{\gamma}-\rho_l^{\gamma})}{\gamma\big(\rho_l\rho_r-\epsilon_1(\rho_l+\rho_r)\big)}}.
\end{align}

Using the Lax entropy inequalities, one can get that the backward shock
wave satisfies
\begin{align}\label{eq4.14}
\sigma^{\epsilon_1\epsilon_2}<\lambda_1(\rho_l,u_l), \quad
\lambda_1(\rho_r,u_r)<\sigma^{\epsilon_1\epsilon_2}<\lambda_2(\rho_r,u_r),
\end{align}
and the forward shock wave satisfies
\begin{align}\label{eq4.15}
\lambda_1(\rho_l,u_l)<\sigma^{\epsilon_1\epsilon_2}<\lambda_2(\rho_l,u_l),
\quad \lambda_2(\rho_r,u_r)<\sigma^{\epsilon_1\epsilon_2}.
\end{align}
Then we can obtain that the following inequality holds for the
backward shock wave
\begin{align}\label{eq4.16}
\dis\frac{-\sqrt{\epsilon_2\rho_r^{\gamma-2}(\rho_r-2\epsilon_1)}}{\rho_l-2\epsilon_1}<\frac{u_r-u_l}{\rho_r-\rho_l}<
\frac{-\sqrt{\epsilon_2\rho_l^{\gamma-2}(\rho_l-2\epsilon_1)}}{\rho_r-2\epsilon_1},
\end{align}
which implies that $\rho_l<\rho_r$ and $u_r<u_l$.

In a analogous way, for the forward shock wave, we have
\begin{align}\label{eq4.17}
\dis\frac{\sqrt{\epsilon_2\rho_r^{\gamma-2}(\rho_r-2\epsilon_1)}}{\rho_l-2\epsilon_1}<\frac{u_r-u_l}{\rho_r-\rho_l}<
\frac{\sqrt{\epsilon_2\rho_l^{\gamma-2}(\rho_l-2\epsilon_1)}}{\rho_r-2\epsilon_1},
\end{align}
which gives $\rho_l>\rho_r$ and $u_r<u_l$.

Thus, given a left state $(\rho_-,u_-)$, one can get the backward shock wave curve
\begin{align}\label{eq4.18}
\overleftarrow{S}(\rho_-,u_-):\ \
u-u_-=\dis-\sqrt{\frac{\epsilon_2(\rho-\rho_-)(\rho^{\gamma}-\rho_-^{\gamma})}{\gamma\big(\rho_-\rho-\epsilon_1(\rho+\rho_-)\big)}},\
\ \ \rho>\rho_-,
\end{align}
and the forward shock wave curve
\begin{align}\label{eq4.19}
\overrightarrow{S}(\rho_-,u_-):\ \
u-u_-=\dis-\sqrt{\frac{\epsilon_2(\rho-\rho_-)(\rho^{\gamma}-\rho_-^{\gamma})}{\gamma\big(\rho_-\rho-\epsilon_1(\rho+\rho_-)\big)}},\
\ \ \rho<\rho_-.
\end{align}

In addition, for the backward shock wave, differentiating $u$ with
respect to $\rho$ in \eqref{eq4.18}, it is immediate that
\begin{align}\label{eq4.20}
u_{\rho}=-\frac{\epsilon_2}{2}\Bigg(\frac{\epsilon_2(\rho-\rho_-)(\rho^{\gamma}-\rho_-^{\gamma})}{\gamma\big(\rho_-\rho-\epsilon_1(\rho+\rho_-)\big)}\Bigg)^{-\frac{1}{2}}
\frac{I}{\gamma\big(\rho_-\rho-\epsilon_1(\rho+\rho_-)\big)^2}<0,
\end{align}
where
$I=\rho_-(\rho_--2\epsilon_1)(\rho^{\gamma}-\rho_-^{\gamma})+\gamma\rho^{\gamma-1}\big(\rho_-\rho-\epsilon_1(\rho+\rho_-)\big)(\rho-\rho_-)$.
Similarly, for the forward shock wave, we have $u_{\rho}>0$.

When $\rho\rightarrow +\infty$ in \eqref{eq4.18}, we find
$\underset{\rho\rightarrow +\infty}\lim u =-\infty$. When
$\rho\rightarrow 2\epsilon_1$ in \eqref{eq4.19}, we obtain
\begin{align}\label{eq4.21}
\lim_{\rho\rightarrow
2\epsilon_1}u=u_--\dis\sqrt{\frac{\epsilon_2(\rho_-^{\gamma}-(2\epsilon_1)^{\gamma})}{\gamma\epsilon_1}},
\end{align}
which shows that the forward shock wave curve intersects with the
line $\rho=2\epsilon_1$ at the point
$(2\epsilon_1,u_2)=(2\epsilon_1,u_--\sqrt{\frac{\epsilon_2(\rho_-^{\gamma}-(2\epsilon_1)^{\gamma})}{\gamma\epsilon_1}})$.

Through the analysis above, as illustrated in Fig. 3, fixing a left state $(\rho_-,u_-)$, the
phase plane can be divided into five regions by the wave curves.
\begin{center}
%\ resizebox {30}{20}{\includegraphics{fig2.2.1-2.2.2.eps}}
\includegraphics*[134,492][408,639]{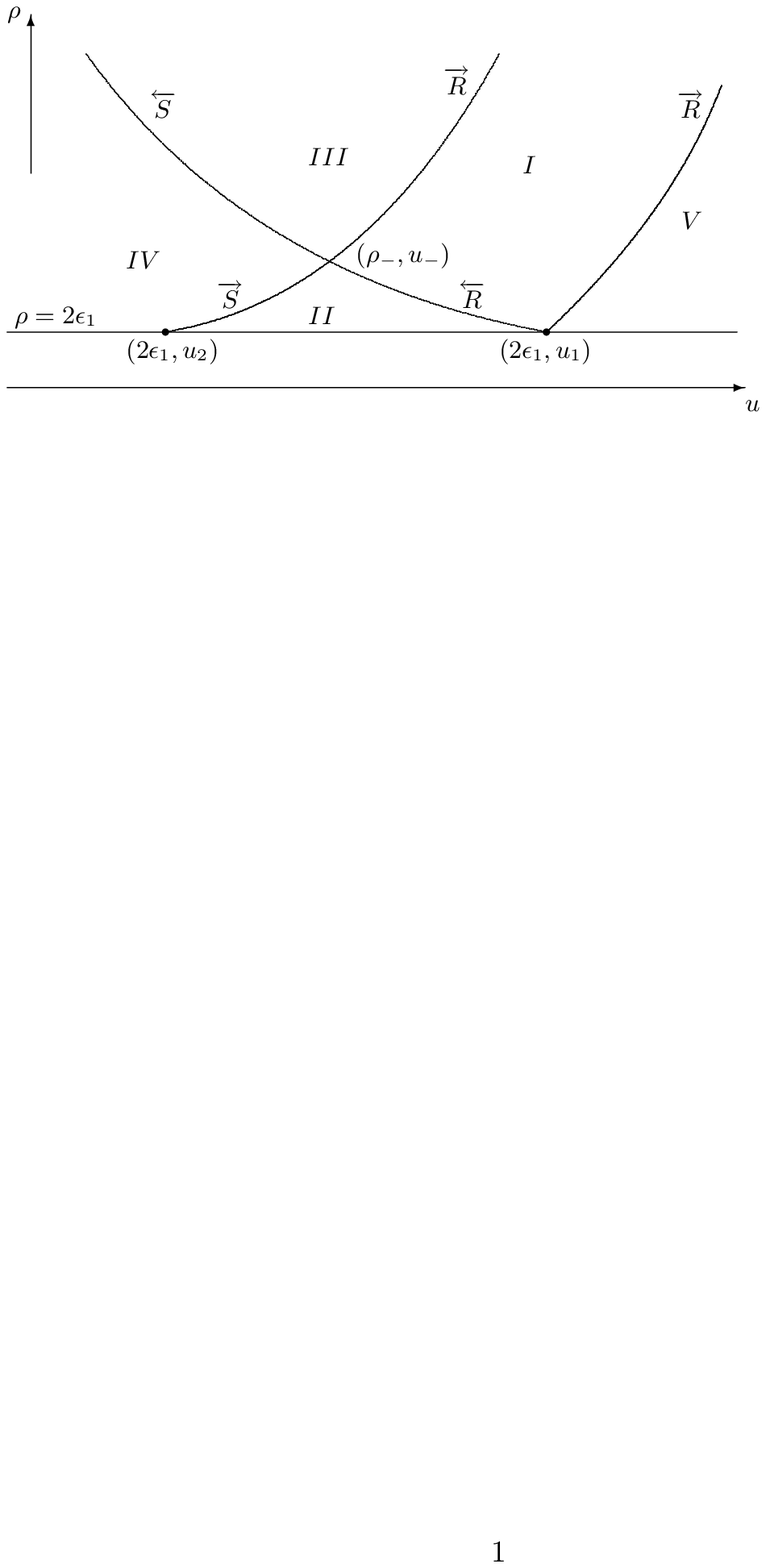}

{\text{ Fig. 3.} Curves of elementary waves.}

\end{center}

Now, according to the right state $(\rho_+,u_+)$ in the different
regions, one can get five kinds of configurations of solutions.
Particularly, when $(\rho_+,u_+)\in
\overleftarrow{S}\overrightarrow{S}(\rho_-,u_-)$, the Riemann
solution contains two shock waves and a nonvacuum intermediate
constant states whose density may become singular as
$\epsilon_1,\epsilon_2\rightarrow0$. When $(\rho_+,u_+)\in
\overleftarrow{R}\overrightarrow{R}(\rho_-,u_-)$, the Riemann
solution contains two rarefaction waves and a  intermediate state
that may be a constant density solution $(\rho=2\epsilon_1)$. Since
the other two regions
$\overleftarrow{S}\overrightarrow{R}(\rho_-,u_-)$ and
$\overleftarrow{R}\overrightarrow{S}(\rho_-,u_-)$ have empty
interiors when $\epsilon_1,\epsilon_2\rightarrow0$, it suffices to
study the limit process for the two cases $(\rho_+,u_+)\in
\overleftarrow{S}\overrightarrow{S}(\rho_-,u_-)$ and
$(\rho_+,u_+)\in \overleftarrow{R}\overrightarrow{R}(\rho_-,u_-)$.

\section {\bf Formation of delta shock waves for the system (1.5)}

This section analyzes the limit as
$\epsilon_1,\epsilon_2\rightarrow0$ of solutions of
\eqref{eq1.5} and \eqref{eq1.7} in the case  $(\rho_+,u_+)\in
\overleftarrow{S}\overrightarrow{S}(\rho_-,u_-)$ with $u_->u_+$.

\noindent \textbf{5.1. Limit behavior of the Riemann solutions as
$\epsilon_1,\epsilon_2\rightarrow0$}

For any $\epsilon_1,\epsilon_2>0$, let
$(\rho_*^{\epsilon_1\epsilon_2},u_*^{\epsilon_1\epsilon_2})$ be the
intermediate state in the sense that $(\rho_-,u_-)$ and
$(\rho_*^{\epsilon_1\epsilon_2},u_*^{\epsilon_1\epsilon_2})$ are
connected by backward shock wave $\overleftarrow{S}$ with speed
$\sigma_1^{\epsilon_1\epsilon_2}$ and that
$(\rho_*^{\epsilon_1\epsilon_2},u_*^{\epsilon_1\epsilon_2})$ and
$(\rho_+,u_+)$ are connected by forward shock wave
$\overrightarrow{S}$ with speed $\sigma_2^{\epsilon_1\epsilon_2}$.
They have the following relations
\begin{align}\label{eq5.1}
\dis
u_*^{\epsilon_1\epsilon_2}-u_-=-\sqrt{\frac{\epsilon_2(\rho_*^{\epsilon_1\epsilon_2}-\rho_-)((\rho_*^{\epsilon_1\epsilon_2})^{\gamma}-\rho_-^\gamma)}{\gamma\big(\rho_-\rho_*^{\epsilon_1\epsilon_2}-\epsilon_1(\rho_-+\rho_*^{\epsilon_1\epsilon_2})\big)}},\
\ \ \ \ \ \ \rho_*^{\epsilon_1\epsilon_2}>\rho_-
\end{align}
on $\overleftarrow{S}$, and
\begin{align}\label{eq5.2}
\dis
u_+-u_*^{\epsilon_1\epsilon_2}=-\sqrt{\frac{\epsilon_2(\rho_+-\rho_*^{\epsilon_1\epsilon_2})(\rho_+^{\gamma}-(\rho_*^{\epsilon_1\epsilon_2})^\gamma)}{\gamma\big(\rho_*^{\epsilon_1\epsilon_2}\rho_+-\epsilon_1(\rho_*^{\epsilon_1\epsilon_2}+\rho_+)\big)}},\
\ \ \ \ \ \ \rho_*^{\epsilon_1\epsilon_2}>\rho_+
\end{align}
on $\overrightarrow{S}$. Then we have the following lemmas.

\vspace{0.2cm}

 \noindent{\textbf{Lemma 5.1.}}\textit{
$\lim\limits_{\epsilon_1,\epsilon_2\rightarrow0}\rho_*^{\epsilon_1\epsilon_2}=+\infty$.}

\vspace{0.2cm}

\textbf{Proof. } Suppose that
$\lim\limits_{\epsilon_1,\epsilon_2\rightarrow0}\rho_*^{\epsilon_1\epsilon_2}=M\in(\text{max}(\rho_-,\rho_+),+\infty)$.
It follows from \eqref{eq5.1} and \eqref{eq5.2} that
\begin{align}\label{eq5.3}
\dis
u_+-u_-=-\sqrt{\frac{\epsilon_2}{\gamma}}\Bigg(\sqrt{\frac{(\rho_*^{\epsilon_1\epsilon_2}-\rho_-)((\rho_*^{\epsilon_1\epsilon_2})^{\gamma}-\rho_-^\gamma)}{\rho_-\rho_*^{\epsilon_1\epsilon_2}-\epsilon_1(\rho_-+\rho_*^{\epsilon_1\epsilon_2})}}
+\sqrt{\frac{(\rho_+-\rho_*^{\epsilon_1\epsilon_2})(\rho_+^{\gamma}-(\rho_*^{\epsilon_1\epsilon_2})^\gamma)}{\rho_*^{\epsilon_1\epsilon_2}\rho_+-\epsilon_1(\rho_*^{\epsilon_1\epsilon_2}+\rho_+)}}\Bigg).
\end{align}
Letting $\epsilon_1,\epsilon_2\rightarrow0$ in \eqref{eq5.3}, one
can get $u_+=u_-$, which contradicts  $u_+<u_-$. Therefore, Lemma
5.1 holds. \hspace{2cm} \ding {122}

\vspace{0.2cm}

Letting $\epsilon_1,\epsilon_2\rightarrow0$ in
\eqref{eq5.3}, one can directly get

\noindent {\textbf{Lemma 5.2.}}
\textit{$\dis\lim\limits_{\epsilon_1,\epsilon_2\rightarrow0}\epsilon_2(\rho_*^{\epsilon_1\epsilon_2})^{\gamma}=\gamma\rho_-\rho_+\Bigg(\frac{u_--u_+}{\sqrt{\rho_-}+\sqrt{\rho_+}}\Bigg)^2.$}

\vspace{0.2cm}

\noindent {\textbf{Lemma 5.3.}} Set $\sigma=\dis
\frac{\sqrt{\rho_-}u_-+\sqrt{\rho_+}u_+}{\sqrt{\rho_-}+\sqrt{\rho_+}}$.
Then \textit{
\begin{align}\label{eq5.4}\lim\limits_{\epsilon_1,\epsilon_2\rightarrow0}u_*^{\epsilon_1\epsilon_2}=\lim\limits_{\epsilon_1,\epsilon_2\rightarrow0}\sigma_1^{\epsilon_1\epsilon_2}=\lim\limits_{\epsilon_1,\epsilon_2\rightarrow0}\sigma_2^{\epsilon_1\epsilon_2}=
\sigma.\end{align}}

\vspace{0.2cm}

\textbf{Proof.} Passing to the limit
$\epsilon_1,\epsilon_2\rightarrow0$ in \eqref{eq5.1} and noticing
Lemma 5.2, we have
\begin{align}\label{eq5.5}
\lim\limits_{\epsilon_1,\epsilon_2\rightarrow0}u_*^{\epsilon_1\epsilon_2}=u_--\frac{1}{\sqrt{\gamma
\rho_-}}\dis\lim\limits_{\epsilon_1,\epsilon_2\rightarrow0}
\sqrt{\epsilon_2(\rho_*^{\epsilon_1\epsilon_2})^{\gamma}}=\sigma.
\end{align}

Form \eqref{eq4.11}, $\sigma_1^{\epsilon_1\epsilon_2}$ and
$\sigma_2^{\epsilon_1\epsilon_2}$ can be calculated by
\begin{align}\label{eq5.6}
\dis\sigma_1^{\epsilon_1\epsilon_2}=u_*^{\epsilon_1\epsilon_2}+\frac{(\rho_--2\epsilon_1)(u_*^{\epsilon_1\epsilon_2}-u_-)}{\rho_*^{\epsilon_1\epsilon_2}-\rho_-},\quad
\dis\sigma_2^{\epsilon_1\epsilon_2}=u_*^{\epsilon_1\epsilon_2}+\frac{(\rho_+-2\epsilon_1)(u_+-u_*^{\epsilon_1\epsilon_2})}{\rho_+-\rho_*^{\epsilon_1\epsilon_2}},
\end{align}
thus,
$\lim\limits_{\epsilon_1,\epsilon_2\rightarrow0}\sigma_1^{\epsilon_1\epsilon_2}=\lim\limits_{\epsilon_1,\epsilon_2\rightarrow0}\sigma_2^{\epsilon_1\epsilon_2}=\lim\limits_{\epsilon_1,\epsilon_2\rightarrow0}u_*^{\epsilon_1\epsilon_2}.$
So the lemma is true. \hspace{2cm} \ding {122}

\vspace{0.2cm}

Lemma 5.1 and Lemma 5.3 show that when $\epsilon_1$ and $\epsilon_2$
drop to zero,  $\overleftarrow{S}$ and $\overrightarrow{S}$
coincide,  the intermediate density $\rho_*^{\epsilon_1\epsilon_2}$
becomes singular.

\vspace{0.2cm}

Combining \eqref{eq5.6} with Lemma 5.1 and Lemma 5.3, we have the
following result.

\vspace{0.2cm}

\noindent {\textbf{Lemma 5.4.}} \textit{$
\lim\limits_{\epsilon_1,\epsilon_2\rightarrow0}\rho_*^{\epsilon_1\epsilon_2}(\sigma_2^{\epsilon_1\epsilon_2}-\sigma_1^{\epsilon_1\epsilon_2})=\sigma[\rho]-[\rho
u]. $}

\vspace{0.2cm}

\noindent \textbf{5.2. Weighted delta shock waves}

Now, we  show the theorem characterizing the limit as
$\epsilon_1,\epsilon_2\rightarrow0$ for the case $u_+< u_-$ and
$(\rho_+,u_+)\in \overleftarrow{S}\overrightarrow{S}(\rho_-,u_-)$.

\vspace{0.2cm}

\noindent {\textbf{Theorem 5.5.}} \textit{Let  $u_+< u_-$.  Assume
$(\rho^{\epsilon_1\epsilon_2},u^{\epsilon_1\epsilon_2})$ is a
two-shock wave solution of $\eqref{eq1.5}$ and $\eqref{eq1.7}$
constructed in Section 4. Then, when
$\epsilon_1,\epsilon_2\rightarrow0$, $\rho^{\epsilon_1\epsilon_2}$
and $\rho^{\epsilon_1\epsilon_2}u^{\epsilon_1\epsilon_2}$ converge
in the sense of distributions, and the limit functions of
$\rho^{\epsilon_1\epsilon_2}$ and $\rho^{\epsilon_1\epsilon_2}
u^{\epsilon_1\epsilon_2}$ are the sum of a step function and a
$\delta$-function with the weights \
$$\dis\frac{t}{\sqrt{1+\sigma^2}}(\sigma[\rho]-[\rho u])\ \ and \ \ \dis\frac{t}{\sqrt{1+\sigma^2}}(\sigma[\rho u]-[\rho u^2]),$$
respectively, which form a delta shock solution of $\eqref{eq1.1}$
with the
 Riemann data $\eqref{eq1.7}$.}

\vspace{0.2cm}

\textbf{Proof.} (i). Set $\xi=x/t$. Then, for each
$\epsilon_1,\epsilon_2>0$, the
 Riemann solution containing two shocks can be expressed as
\begin{align}\label{eq5.7}
(\rho^{\epsilon_1\epsilon_2},u^{\epsilon_1\epsilon_2})(\xi)=\left\{\begin{array}{ll}
 (\rho_-,u_-),&\xi<\sigma_1^{\epsilon_1\epsilon_2}, \\[2mm]
 (\rho_*^{\epsilon_1\epsilon_2},u_*^{\epsilon_1\epsilon_2}), &\sigma_1^{\epsilon_1\epsilon_2}<\xi<\sigma_2^{\epsilon_1\epsilon_2},\\[2mm]
 (\rho_+,u_+), &\xi>\sigma_2^{\epsilon_1\epsilon_2},
 \end{array}\right.
\end{align}
satisfying weak formulations: For any $\phi\in
C^1_0(-\infty,+\infty)$,
\begin{align}\label{eq5.8}
\dis\int^{+\infty}_{-\infty}(\rho^{\epsilon_1\epsilon_2}
u^{\epsilon_1\epsilon_2}-\rho^{\epsilon_1\epsilon_2} \xi-2\epsilon_1
u^{\epsilon_1\epsilon_2})\phi'd\xi-\dis\int^{+\infty}_{-\infty}\rho^{\epsilon_1\epsilon_2}\phi
d\xi=0,
\end{align}
and
\begin{align}\label{eq5.9}
\dis\int^{+\infty}_{-\infty}\Big((\rho^{\epsilon_1\epsilon_2}-\epsilon_1)
(u^{\epsilon_1\epsilon_2})^2+\frac{\epsilon_2(\rho^{\epsilon_1\epsilon_2})^{\gamma}}{\gamma}-\rho^{\epsilon_1\epsilon_2}
u^{\epsilon_1\epsilon_2}
\xi\Big)\phi'd\xi-\dis\int^{+\infty}_{-\infty}\rho^{\epsilon_1\epsilon_2}
u^{\epsilon_1\epsilon_2} \phi d\xi=0.
\end{align}

(ii). The first integral in \eqref{eq5.8} can be decomposed into
\begin{align}\label{eq5.10}
\bigg(\dis\int^{\sigma_1^{\epsilon_1\epsilon_2}}_{-\infty}+\dis\int_{\sigma_1^{\epsilon_1\epsilon_2}}^{\sigma_2^{\epsilon_1\epsilon_2}}
+\dis\int^{+\infty}_{\sigma_2^{\epsilon_1\epsilon_2}}\bigg)(\rho^{\epsilon_1\epsilon_2}
u^{\epsilon_1\epsilon_2}-\rho^{\epsilon_1\epsilon_2} \xi-2\epsilon_1
u^{\epsilon_1\epsilon_2})\phi'd\xi.
\end{align}
The limit of the sum of the first and last term of \eqref{eq5.10} equals
\begin{align}\label{eq5.11}
\begin{array}{l}
\lim\limits_{\epsilon_1,\epsilon_2\rightarrow0}\dis\int^{\sigma_1^{\epsilon_1\epsilon_2}}_{-\infty}(\rho_- u_--\rho_- \xi-2\epsilon_1
u_-)\phi'd\xi+
 \lim\limits_{\epsilon_1,\epsilon_2\rightarrow0}\dis\int_{\sigma_2^{\epsilon_1\epsilon_2}}^{+\infty}(\rho_+ u_+-\rho_+ \xi-2\epsilon_1
u_+)\phi'd\xi\\[4mm]
 =(\sigma[\rho]-[\rho u])\phi(\sigma)+\dis\int^{+\infty}_{-\infty}H(\xi-\sigma)\phi
 d\xi
\end{array}
\end{align}
with $H(\xi-\sigma)$ taking $\rho_-$ for $\xi<\sigma$ and $\rho_+$ for $\xi>\sigma$, respectively.
While the limit of the second term  of \eqref{eq5.10} can be written as
\begin{align}\label{eq5.12}
\begin{array}{ll}
&\lim\limits_{\epsilon_1,\epsilon_2\rightarrow0}\dis\int_{\sigma_1^{\epsilon_1\epsilon_2}}^{\sigma_2^{\epsilon_1\epsilon_2}}(\rho_*^{\epsilon_1\epsilon_2} u_*^{\epsilon_1\epsilon_2}-\rho_*^{\epsilon_1\epsilon_2} \xi-2\epsilon_1 u_*^{\epsilon_1\epsilon_2})\phi'd\xi\\[5mm]
  &=\lim\limits_{\epsilon_1,\epsilon_2\rightarrow0}\rho_*^{\epsilon_1\epsilon_2}(\sigma_2^{\epsilon_1\epsilon_2}-\sigma_1^{\epsilon_1\epsilon_2})
     \bigg(\frac{\phi(\sigma_2^{\epsilon_1\epsilon_2})-\phi(\sigma_1^{\epsilon_1\epsilon_2})}{\sigma_2^{\epsilon_1\epsilon_2}-\sigma_1^{\epsilon_1\epsilon_2}}u_*^{\epsilon_1\epsilon_2}
 -\frac{\sigma_2^{\epsilon_1\epsilon_2}\phi(\sigma_2^{\epsilon_1\epsilon_2})-\sigma_1^{\epsilon_1\epsilon_2}\phi(\sigma_1^{\epsilon_1\epsilon_2})}{\sigma_2^{\epsilon_1\epsilon_2}-\sigma_1^{\epsilon_1\epsilon_2}}\\[5mm]
&\hspace{2.0cm}+\frac{1}{\sigma_2^{\epsilon_1\epsilon_2}-\sigma_1^{\epsilon_1\epsilon_2}}\dis\int_{\sigma_1^{\epsilon_1\epsilon_2}}^{\sigma_2^{\epsilon_1\epsilon_2}}\phi
     d\xi\bigg)
-\lim\limits_{\epsilon_1,\epsilon_2\rightarrow0}2\epsilon_1u_*^{\epsilon_1\epsilon_2}\Big(\phi(\sigma_2^{\epsilon_1\epsilon_2})-\phi(\sigma_1^{\epsilon_1\epsilon_2})\Big)\\[5mm]
  &=(\sigma[\rho]-[\rho u])\Big(\sigma\phi'(\sigma)-\sigma\phi'(\sigma)-\phi(\sigma)+\phi(\sigma)\Big)\\[3mm]
  & =0.
\end{array}
\end{align}
Returning to \eqref{eq5.8}, we immediately obtain that
\begin{align}\label{eq5.13}
 \lim\limits_{\epsilon_1,\epsilon_2\rightarrow0}\dis\int^{+\infty}_{-\infty}\rho^{\epsilon_1\epsilon_2}\phi d\xi
 =(\sigma[\rho]-[\rho u])\phi(\sigma)+\dis\int^{+\infty}_{-\infty}H(\xi-\sigma)\phi
 d\xi.
\end{align}

(iii) Now we consider the limit of $\rho^{\epsilon_1\epsilon_2}
u^{\epsilon_1\epsilon_2}$. In the same way as before, we decompose the first
integral of \eqref{eq5.9} into
\begin{align}\label{eq5.14}
\bigg(\dis\int^{\sigma_1^{\epsilon_1\epsilon_2}}_{-\infty}+\dis\int_{\sigma_1^{\epsilon_1\epsilon_2}}^{\sigma_2^{\epsilon_1\epsilon_2}}
+\dis\int^{+\infty}_{\sigma_2^{\epsilon_1\epsilon_2}}\bigg)\Big((\rho^{\epsilon_1\epsilon_2}-\epsilon_1)
(u^{\epsilon_1\epsilon_2})^2+\frac{\epsilon_2(\rho^{\epsilon_1\epsilon_2})^{\gamma}}{\gamma}-\rho^{\epsilon_1\epsilon_2}
u^{\epsilon_1\epsilon_2} \xi\Big)\phi'd\xi.
\end{align}
As $\epsilon_1,\epsilon_2\rightarrow0$, the limit of the sum of the first and last term of \eqref{eq5.14} is
\begin{align}\label{eq5.15}
(\sigma[\rho u]-[\rho
u^2])\phi(\sigma)+\dis\int^{+\infty}_{-\infty}\widetilde{H}(\xi-\sigma)\phi d\xi
\end{align}
with $\widetilde{H}(\xi-\sigma)$ taking $\rho_-u_-$ for $\xi<\sigma$ and $\rho_+u_+$ for $\xi>\sigma$, respectively.
Applying Lemmas 5.1-5.4, one can deduce that the limit of the second term  of \eqref{eq5.14} equals
\begin{align}\label{eq5.16}
\begin{array}{ll}
&\lim\limits_{\epsilon_1,\epsilon_2\rightarrow0}\rho_*^{\epsilon_1\epsilon_2}(\sigma_2^{\epsilon_1\epsilon_2}-\sigma_1^{\epsilon_1\epsilon_2})
\bigg(\textstyle
\frac{\phi(\sigma_2^{\epsilon_1\epsilon_2})-\phi(\sigma_1^{\epsilon_1\epsilon_2})}{\sigma_2^{\epsilon_1\epsilon_2}-\sigma_1^{\epsilon_1\epsilon_2}}(u_*^{\epsilon_1\epsilon_2})^2
+\frac{\epsilon_2(\rho_*^{\epsilon_1\epsilon_2})^{\gamma-1}}{\gamma}\frac{\phi(\sigma_2^{\epsilon_1\epsilon_2})-\phi(\sigma_1^{\epsilon_1\epsilon_2})}{\sigma_2^{\epsilon_1\epsilon_2}
-\sigma_1^{\epsilon_1\epsilon_2}}\\[5mm]
&\hspace{3cm}-\frac{\sigma_2^{\epsilon_1\epsilon_2}\phi(\sigma_2^{\epsilon_1\epsilon_2})-\sigma_1^{\epsilon_1\epsilon_2}\phi(\sigma_1^{\epsilon_1\epsilon_2})}{\sigma_2^{\epsilon_1\epsilon_2}-\sigma_1^{\epsilon_1\epsilon_2}}u_*^{\epsilon_1\epsilon_2}
+\frac{u_*^{\epsilon_1\epsilon_2}}{\sigma_2^{\epsilon_1\epsilon_2}-\sigma_1^{\epsilon_1\epsilon_2}}\dis\int_{\sigma_1^{\epsilon_1\epsilon_2}}^{\sigma_2^{\epsilon_1\epsilon_2}}\phi
d\xi\bigg)\\[5mm]
&\hspace{1.0cm}-\lim\limits_{\epsilon_1,\epsilon_2\rightarrow0}\epsilon_1(u_*^{\epsilon_1\epsilon_2})^2\Big(\phi(\sigma_2^{\epsilon_1\epsilon_2})-\phi(\sigma_1^{\epsilon_1\epsilon_2})\Big)=0.
\end{array}
\end{align}
Thus, it follows from \eqref{eq5.9} that
\begin{align}\label{eq5.17}
 \lim\limits_{\epsilon_1,\epsilon_2\rightarrow0}\dis\int^{+\infty}_{-\infty}\rho^{\epsilon_1\epsilon_2} u^{\epsilon_1\epsilon_2}\phi d\xi
=(\sigma[\rho u]-[\rho
u^2])\phi(\sigma)+\dis\int^{+\infty}_{-\infty}\widetilde{H}(\xi-\sigma)\phi
d\xi.
\end{align}

(iiii). Finally, we analyze the limit of
$\rho^{\epsilon_1\epsilon_2}$ and
$\rho^{\epsilon_1\epsilon_2}u^{\epsilon_1\epsilon_2}$ by tracking
the time-dependence of the weights of the $\delta$-measures as
$\epsilon_1,\epsilon_2\rightarrow0$.

Taking \eqref{eq5.13} into account, we have for any $\psi\in
C^\infty_0(R\times R^+)$
\begin{align}\label{eq5.18}
\begin{array}{l}
\lim\limits_{\epsilon_1,\epsilon_2\rightarrow0}\dis\int^{+\infty}_{0}\dis\int^{+\infty}_{-\infty}\rho^{\epsilon_1\epsilon_2}(x/t)\psi(x,t)dxdt\\[5mm]
 =\lim\limits_{\epsilon_1,\epsilon_2\rightarrow0}\dis\int^{+\infty}_{0}t\bigg(\dis\int^{+\infty}_{-\infty}\rho^{\epsilon_1\epsilon_2}(\xi)\psi(\xi t,t)d\xi\bigg)dt\\[5mm]
 =\dis\int^{+\infty}_{0}(\sigma[\rho]-[\rho u])t\psi(\sigma t,t)dt+\dis\int^{+\infty}_{0}\dis\int^{+\infty}_{-\infty}H(x-\sigma
 t)\psi(x,t)dxdt,
\end{array}
\end{align}
in which, by the definition \eqref{eq2.2}, we get
\begin{align}\label{eq5.19}
 \dis\int^{+\infty}_{0}(\sigma[\rho]-[\rho u])t\psi(\sigma t,t)dt=\Big<w_1(\cdot)\delta_S,\psi(\cdot,\cdot)\Big>
\end{align}
with $w_1(t)=\dis\frac{t}{\sqrt{1+\sigma^2}}(\sigma[\rho]-[\rho u])$. Similarly, one can show that
\begin{align}\label{eq5.20}
\begin{array}{l}
\lim\limits_{\epsilon_1,\epsilon_2\rightarrow0}\dis\int^{+\infty}_{0}\dis\int^{+\infty}_{-\infty}\rho^{\epsilon_1\epsilon_2}u^{\epsilon_1\epsilon_2}(x/t)\psi(x,t)dxdt\\[5mm]
\qquad =\Big<w_2(\cdot)\delta_S,\psi(\cdot,\cdot)\Big>+\dis\int^{+\infty}_{0}\dis\int^{+\infty}_{-\infty}\widetilde{H}(x-\sigma
 t)\psi(x,t)dxdt
\end{array}
\end{align}
with $w_2(t)=\dis\frac{t}{\sqrt{1+\sigma^2}}(\sigma[\rho u]-[\rho u^2])$.

\noindent The proof of Theorem 5.5 is finished.  \hspace{2cm} \ding {122}

\section {\bf  Formation of vacuum states for the system (1.5)}

In this section, we study the limit of the solutions of \eqref{eq1.5} and
\eqref{eq1.7} in the case $(\rho_+,u_+)\in
\overleftarrow{R}\overrightarrow{R}(\rho_-,u_-)$ with $u_-<u_+$ as
$\epsilon_1,\epsilon_2\rightarrow0$.

\vspace{0.2cm}

According to Section 4, one can get that, on the backward centred
rarefaction wave, the solution satisfies
\begin{align}\label{eq6.1}
\left\{\begin{array}{l}
\xi=u^{\epsilon_1\epsilon_2}-\sqrt{\epsilon_2
(\rho^{\epsilon_1\epsilon_2})^{\gamma-2}(\rho^{\epsilon_1\epsilon_2}-2\epsilon_1)},\cr\noalign
{\vskip2truemm} u_--\sqrt{\epsilon_2
\rho_-^{\gamma-2}(\rho_--2\epsilon_1)}<\xi<u_*^{\epsilon_1\epsilon_2}-\sqrt{\epsilon_2
(\rho_*^{\epsilon_1\epsilon_2})^{\gamma-2}(\rho_*^{\epsilon_1\epsilon_2}-2\epsilon_1)},\quad
\rho_*^{\epsilon_1\epsilon_2}<\rho_-,
\end{array}\right.
\end{align}
and, on the forward centred rarefaction wave,
\begin{align}\label{eq6.2}
\left\{\begin{array}{l}
\xi=u^{\epsilon_1\epsilon_2}+\sqrt{\epsilon_2
(\rho^{\epsilon_1\epsilon_2})^{\gamma-2}(\rho^{\epsilon_1\epsilon_2}-2\epsilon_1)},\cr\noalign
{\vskip2truemm}u_*^{\epsilon_1\epsilon_2}+\sqrt{\epsilon_2
(\rho_*^{\epsilon_1\epsilon_2})^{\gamma-2}(\rho_*^{\epsilon_1\epsilon_2}-2\epsilon_1)}<\xi<u_++\sqrt{\epsilon_2
\rho_+^{\gamma-2}(\rho_+-2\epsilon_1)},\quad
\rho_*^{\epsilon_1\epsilon_2}<\rho_+.
\end{array}\right.
\end{align}

Now, we can conclude the following theorem.

\vspace{0.2cm}

\noindent {\textbf{Theorem 6.1.}} \textit{Let  $u_-<u_+$.  Assume
$(\rho^{\epsilon_1\epsilon_2},u^{\epsilon_1\epsilon_2})$ is a
two-rarefaction wave solution of $\eqref{eq1.5}$ and $\eqref{eq1.7}$
constructed in Section 4. Then, there exist $\epsilon_0>0$, when $0<\epsilon_1<\epsilon_0$ and
$0<\epsilon_2<\epsilon_0$,  the constant density solution
$(\rho=2\epsilon_1)$ appears in the solution. And as
$\epsilon_1,\epsilon_2 \rightarrow0$, the two rarefaction waves become two contact
discontinuities connecting the constant states $(u_\pm,\rho_\pm)$
and the vacuum ($\rho=0$), which form a vacuum solution of \eqref{eq1.1} with the Riemann data $\eqref{eq1.7}$.}

\vspace{0.2cm}

\textbf{Proof.} Set $\epsilon_1=\epsilon_2=\epsilon_0$. Since
$(\rho_*^{\epsilon_1\epsilon_2},u_*^{\epsilon_1\epsilon_2})$ is on
the curve $\overleftarrow{R}(\rho_-,u_-)$, we have
\begin{align}\label{eq6.3}
u_*^{\epsilon_1\epsilon_2}=u_--\dis\int^{\rho_*^{\epsilon_1\epsilon_2}}_{\rho_-}
\sqrt{\frac{\epsilon_0 s^{\gamma-2}}{s-2\epsilon_0}} ds \leq
u_-+\dis\int^{\rho_-}_{2\epsilon_0} \sqrt{\frac{\epsilon_0
s^{\gamma-2}}{s-2\epsilon_0}} ds=A^{\epsilon_0}.
\end{align}

When $u_-<u_+<A^{\epsilon_0}$, there is no constant-density
 in the solution. That is,
there exist $\epsilon_{01}$ such that $(\rho_+,u_+)\in
I(\rho_-,u_-)$ when $u_-<u_+<A^{\epsilon_{01}}$.

However, when $A^{\epsilon_0}<u_+$, the constant density solution
appears, which implies that there exist $\epsilon_{02}$ such that
$(\rho_+,u_+)\in V(\rho_-,u_-)$ when $A^{\epsilon_{02}}<u_+$.

Let $f(\epsilon)=\dis\int^{\rho_-}_{2\epsilon} \sqrt{\frac{\epsilon
s^{\gamma-2}}{s-2\epsilon}}ds-u_++u_-$. Since the integral $\dis\int^{\rho_-}_{0}\sqrt{\frac{
s^{\gamma-2}+2}{s}}ds$ is convergent, one can deduce that, thanks to M-criterion, the integral $\dis\int^{\rho_-}_{2\epsilon} \sqrt{\frac{\epsilon s^{\gamma-2}}{s-2\epsilon}}ds$
is uniformly convergent in $\varepsilon$, then the function $f(\epsilon)$ is
continuous with respect to $\epsilon$ and
$f(\epsilon_{01})f(\epsilon_{02})<0$. Thus, there exists
$\epsilon_0\in[\epsilon_{02},\epsilon_{01}]$ such that
$f(\epsilon_0)=0$.

So when $0<\epsilon_1<\epsilon_0$ and $0<\epsilon_2<\epsilon_0$, the
density of the intermediate state becomes a constant  with
\begin{align}\label{eq6.4}
(\rho_*^{\epsilon_1\epsilon_2},u_*^{\epsilon_1\epsilon_2})(\xi)=(2\epsilon_1,\xi),
\quad u_1^{\epsilon_1\epsilon_2}\leq \xi \leq
u_2^{\epsilon_1\epsilon_2},
\end{align}
where
$$
u_1^{\epsilon_1\epsilon_2}=u_-+\dis\int^{\rho_-}_{2\epsilon_1}
\sqrt{\frac{\epsilon_2 s^{\gamma-2}}{s-2\epsilon_1}}ds,\quad
u_2^{\epsilon_1\epsilon_2}=u_+-\dis\int^{\rho_+}_{2\epsilon_1}
\sqrt{\frac{\epsilon_2 s^{\gamma-2}}{s-2\epsilon_1}}ds.
$$
Thus, letting $\epsilon_1,\epsilon_2 \rightarrow0$, one can find $\lim\limits_{\epsilon_1,\epsilon_2\rightarrow0}\rho_*^{\epsilon_1\epsilon_2}=0$.
Using the uniform boundedness of $\rho^{\epsilon_1\epsilon_2}$
with respect to $\epsilon_1$ and $\epsilon_2$, it follows that
$$
\lim\limits_{\epsilon_1,\epsilon_2\rightarrow0}u_1^{\epsilon_1\epsilon_2}=u_-,
\quad
\lim\limits_{\epsilon_1,\epsilon_2\rightarrow0}u_2^{\epsilon_1\epsilon_2}=u_+,
$$
$$
\lim\limits_{\epsilon_1,\epsilon_2\rightarrow0}u^{\epsilon_1\epsilon_2}(\xi)=\xi
\quad \text{for}\quad \xi \in (u_-,u_+).
$$

In summary, the limit solution for this case can be expressed as
\eqref{eq2.1}, which is a solution of \eqref{eq1.1} containing two contact discontinuities
$\xi=x/t=u_{\pm}$ and a vacuum state in between. This completes the proof of Theorem 6.1. \hspace{2cm} \ding {122}

\vspace{0.2cm}
\noindent {\textbf{Remark.}} The processes of formation of delta shock waves and vacuum states
can be examined with some numerical results as $\epsilon_1$ and $\epsilon_2$
decrease. The numerical simulations will be presented in the version for publication.

\end{document}